\RequirePackage{ifpdf}
\ifpdf % We are running pdfTeX in pdf mode
\documentclass[pdftex]{sigma}
\else
\documentclass{sigma}
\fi

\numberwithin{equation}{section}
\numberwithin{theorem}{section}
\numberwithin{proposition}{section}
\numberwithin{lemma}{section}
\numberwithin{corollary}{section}
\numberwithin{remark}{section}

\begin{document}

\allowdisplaybreaks

\renewcommand{\thefootnote}{$\star$}

\renewcommand{\PaperNumber}{041}

\FirstPageHeading

\ShortArticleName{A First Order $q$-Dif\/ference System for the $BC_1$-Type Jackson Integral} % and Its Applications}

\ArticleName{A First Order $\boldsymbol{q}$-Dif\/ference System
for the $\boldsymbol{BC_1}$-Type \\
Jackson Integral and Its Applications\footnote{This paper is a contribution to the Proceedings of the Workshop ``Elliptic Integrable Systems, Isomonodromy Problems, and Hypergeometric Functions'' (July 21--25, 2008, MPIM, Bonn, Germany). The full collection
is available at
\href{http://www.emis.de/journals/SIGMA/Elliptic-Integrable-Systems.html}{http://www.emis.de/journals/SIGMA/Elliptic-Integrable-Systems.html}}}

\Author{Masahiko ITO}

\AuthorNameForHeading{M. Ito}

\Address{Department of Physics and Mathematics, Aoyama Gakuin University,\\ Kanagawa 229-8558, Japan}

\Email{\href{mailto:mito@gem.aoyama.ac.jp}{mito@gem.aoyama.ac.jp}}

\ArticleDates{Received December 01, 2008, in f\/inal form March 18,
2009; Published online April 03, 2009}

\Abstract{We present an explicit expression for the $q$-dif\/ference system,
which the $BC_1$-type Jackson integral ($q$-series) satisf\/ies,
as f\/irst order simultaneous $q$-dif\/ference equations with a~concrete basis.
As an application, we give a simple proof for the hypergeometric summation formula introduced by Gustafson
and the product formula of the $q$-integral introduced by Nassrallah--Rahman and Gustafson.}


\Classification{33D15; 33D67; 39A13}

\section{Introduction}

\looseness=1
A lot of summation and transformation formulae
for basic hypergeometric series have been found to date.
The $BC_1$-type Jackson integral, which is the main subject of interest in this paper,
is a $q$-series which can be written as a basic hypergeometric series
in a class of so called {\it very-well-poised-balanced} ${}_{2r}\psi_{2r}$.
A key reason to consider the $BC_1$-type Jackson integrals
is to give an explanation of these hypergeometric series
from the view points of the Weyl group symmetry and the $q$-dif\/ference equations of the $BC_1$-type Jackson integrals.
In \cite{IS}, we showed that
Slater's transformation formula
for a very-well-poised-balanced ${}_{2r}\psi_{2r}$ series
can be regarded as a~connection formula for the solutions of
$q$-dif\/ference equations of the $BC_1$-type Jackson integral,
i.e.,
the Jackson integral as a general solution of
$q$-dif\/ference system is  written as a~linear combination of particular solutions.
As a consequence we gave a simple proof of Slater's transformation formula.
(See \cite{IS} for details.
Also see \cite{Ito9} for a connection formula for the
$BC_n$-type Jackson integral, which is
a multisum generalization of
that of type~$BC_1$.)

The aim of this paper is to present an explicit form of the $q$-dif\/ference system as
f\/irst order simultaneous $q$-dif\/ference equations
for the $BC_1$-type Jackson integral
with generic condition on the parameters.
We give the Gauss decomposition of the coef\/f\/icient matrix of the system with a concrete basis (see Theorem \ref{thm:difference}).
Each entry of the decomposed matrices is written as a product of binomials and,
as a consequence,
the determinant of
the coef\/f\/icient matrix is easy to calculate explicitly.
As an application we give a simple proof of
the product formula for Gustafson's multiple $C_n$-type sum \cite{Gu4}.
We also present an explicit form of the $q$-dif\/ference system
for the $BC_1$-type Jackson integral
with a balancing condition on the parameters.
We f\/inally give a simple proof of
the product formula for the $q$-integral of Nassrallah--Rahman~\cite{NR} and Gustafson~\cite{Gu3}.
A recent work of Rains and Spiridonov~\cite{RS} contains results
for the elliptic hypergeometric integral of a similar type to those contained for the $BC_1$-type Jackson integral obtained here.

\section[$BC_1$-type Jackson integral]{$\boldsymbol{BC_1}$-type Jackson integral}

Throughout this paper, we assume $0<q<1$
and denote the $q$-{\it shifted factorial} for all inte\-gers~$N$ by
$
(x)_\infty :=\prod\limits_{i=0}^\infty (1-q^i x)$
and
$
(x)_N:=(x;q)_\infty /(q^Nx;q)_\infty
$.

Let $\mathcal{O}({\mathbb C}^*)$ be the set of holomorphic functions on
the complex multiplicative group ${\mathbb C}^*$.
A~function $f$ on ${\mathbb C}^*$ is said to be {\it symmetric} or {\it skew-symmetric} under the Weyl group action $z\to z^{-1}$ if
$f$ satisf\/ies $f(z)=f(z^{-1})$ or $f(z)=-f(z^{-1})$, respectively.
For $\xi\in {\mathbb C}^*$ and a function $f$ on ${\mathbb C}^*$,
we def\/ine the sum over the lattice ${\mathbb Z}$
\begin{gather*}
%\label{eq:JI}
\int_0^{\xi\infty}
\!\!\!
f(z)\frac{d_q z}{z}
:=(1-q)\sum_{\nu=-\infty}^{\infty}
f(q^\nu \xi),
\end{gather*}
which, provided the integral converges, we call the {\it Jackson integral}.
For an arbitrary positive integer $s$,
we def\/ine the function $\Phi$
and the skew-symmetric function $\Delta$ on ${\mathbb C}^*$ as follows:
\begin{gather}
\Phi (z):=\prod_{m=1}^{2s+2}
z^{\frac{1}{2}-\alpha _{m}}
\frac{(qz/a_{m})_{\infty }}{(za_{m})_{\infty}},
\qquad
\Delta(z):=z^{-1}-z,
\label{eq:Phi}
\end{gather}
where $a_{m}=q^{\alpha_{m}}$.
For a symmetric function $\varphi$ on ${\mathbb C}^*$ and a point $\xi\in {\mathbb C}^*$,
we def\/ine the following sum over the lattice $\mathbb Z$:
\[
\int_0^{\xi\infty}
\!\!\!
\varphi (z)\Phi(z)\Delta(z)
\frac{d_q z}{z},
\]
which we call
the {\it Jackson integral of type $BC_1$} and is simply denoted by
$\langle
\varphi,\xi
\rangle
$.
By def\/inition
the sum
$
\langle
\varphi,\xi
\rangle
$ is invariant under the shift
$\xi\to q^\nu \xi$
for $\nu\in {\mathbb Z}$.

Let $\Theta(z)$ be the function on ${\mathbb C}^*$ def\/ined by
\[
\Theta(z):=
\frac{z^{s-\alpha _{1}-\cdots -\alpha _{2s+2}}
      \theta(z^2)}
     {\prod\limits_{m=1}^{2s+2}\theta(a_{m}z)},
\]
where $\theta(z)$ denotes the function $(z)_\infty(q/z)_\infty$,
which satisf\/ies
\begin{gather}
\theta(qz)=-\theta(z)/z
\qquad\mbox{and}\qquad
\theta(q/z)=\theta(z).
\label{eq:thProperty}
\end{gather}
For a symmetric function $\varphi\in \mathcal{O}({\mathbb C}^*)$,
we denote the function
$
\langle
\varphi,z
\rangle/\Theta(z)
$
by
$
\langle\!\langle
\varphi,z
\rangle\!\rangle
$.
We call
$
\langle\!\langle
\varphi,z
\rangle\!\rangle
$
the {\it regularized Jackson integral of type} $BC_1$, which satisf\/ies the following:
%%%%%%%%%%%%%%%%%%%%%%%%%%%%%%%%%%%%
\begin{lemma}\label{lem:regular}
%%%%%%%%%%%%%%%%%%%%%%%%%%%%%%%%%%%%
Assume $\alpha_1+\alpha_2+\cdots+\alpha_{2s+2}
\not\in\frac{1}{2}+{\mathbb Z}$.
If $\varphi\in \mathcal{O}({\mathbb C}^*)$ is symmetric,
then the function
$
\langle\!\langle
\varphi,z
\rangle\!\rangle
$
is symmetric and holomorphic on ${\mathbb C}^*$.
\end{lemma}
%%%%%%%%%%%%%%%%%%%%%%%%%%%%%%%%%%%%
%
\begin{proof}
See \cite[Proposition 2.2]{IS}.
\end{proof}
%%%%%%%%%%%%%%%%%%%%%%%%%%%%%%%%%%%%
For  an arbitrary meromorphic function $\varphi$ on ${\mathbb C}^*$
we def\/ine the function $\nabla\varphi$ on ${\mathbb C}^*$ by
\[
\nabla\varphi(z):=\varphi(z)-\frac{\Phi(qz)}{\Phi(z)}\varphi(qz).
\]
In particular, from (\ref{eq:Phi}), the function $\Phi(qz)/\Phi(z)$ is the rational function
\[
\frac{\Phi(qz)}{\Phi(z)}=q^{s+1}\prod_{m=1}^{2s+2}\frac{1-a_mz}{a_m-qz}.
\]
The following proposition will be used for the proof of the key equation (Theorem \ref{thm:key}):
%%%%%%%%%%%%%%%%%%%%%%%%%%%%%%%%%%
\begin{lemma}
\label{lem:cohom=0}
%%%%%%%%%%%%%%%%%%%%%%%%%%%%%%%%%%
If
$\displaystyle
\int_0^{\xi\infty}
\!\!\!
\Phi(z)\varphi (z)
\frac{d_q z}{z}
$
is convergent for $\varphi\in \mathcal{O}({\mathbb C}^*)$, then
$\displaystyle
\int_0^{\xi\infty}
\!\!\!
\Phi(z)\nabla\varphi (z)
\frac{d_q z}{z}=0.
$
\end{lemma}
%%%%%%%%%%%%%%%%%%%%%%%%%%%%%%%%%%
\begin{proof}
 See \cite[Lemma 5.1]{Ito7} for instance.
\end{proof}

\section{Key equation}

In this section, we will present a key equation to construct the dif\/ference equations
for the $BC_1$-type Jackson integral.
Before we state it, we introduce the function $e(x;y)$ def\/ined by
\[e(x;y):=x+x^{-1}-\left(y+y^{-1}\right),\]
which is expressed by the product form
\begin{gather*}
e(x;y)=\frac{(y-x)(1-xy)}{xy}.
\end{gather*}
The basic properties of $e(x;y)$ are the following:
\begin{gather}
\bullet\quad e(x;z)=e(x;y)+e(y;z),
\label{eq:e-1}\\
\bullet\quad e(x;y)=-e(y;x), \quad e(x;y)=e\big(x^{-1};y\big),
\label{eq:e-2}\\
\bullet\quad e(x;y)e(z;w)-e(x;z)e(y;w)+e(x;w)e(y;z)=0.
\label{eq:e-3}
\end{gather}
\begin{remark}
As we will see later, equation~(\ref{eq:e-3}) is ignorable in the case $a_1a_2\cdots a_{2s+2}\ne 1$, while equation~(\ref{eq:e-1}) is ignorable in the case $a_1a_2\cdots a_{2s+2}= 1$.
\end{remark}
For functions $f$, $g$ on ${\mathbb C}^*$, the function $fg$ on ${\mathbb C}^*$ is def\/ined by
\[
(fg)(z):=f(z)g(z)\qquad\mbox{for}\qquad z\in {\mathbb C}^*.
\]
Set $e_i(z):=e(z;a_i)$ and $(e_{i_1}e_{i_2}{\cdots} e_{i_s})(z):=e_{i_1}(z)e_{i_2}(z){\cdots} e_{i_s}(z)$.
The symbol $(e_{i_1}{\cdots} \widehat{e}_{i_k}{\cdots} e_{i_s})(z)$ is equal to $(e_{i_1}\cdots e_{i_{k-1}}e_{i_{k+1}}\cdots e_{i_s})(z)$.
The key equation is the following:
%%%%%%%%%%%%%%%%%%%%%%%%%%%%%%
\begin{theorem}
\label{thm:key}
%%%%%%%%%%%%%%%%%%%%%%%%%%%%%%
Suppose $a_i\ne a_j$ if $i\ne j$. If $\{i_1,i_2,\ldots,i_s\}\subset\{1,2,\ldots,2s+2\}$, then
\begin{gather*}
C_0\langle e_{i_1}e_{i_2}\cdots e_{i_s},\xi\rangle
+\sum_{k=1}^s C_{i_k} \langle e_{i_1}\cdots \widehat{e}_{i_k}\cdots e_{i_s},\xi\rangle=0,
%\label{eq:key}
\end{gather*}
where the coefficients $C_0$ and $C_{i_k}$ $(1\le k\le s)$ are given by
\begin{gather*}
C_0=1-a_1a_2\cdots a_{2s+2}
\qquad\mbox{and}\qquad
C_{i_k}=\frac{\prod\limits_{m=1}^{2s+2}(1-a_{i_k}a_m)}{a_{i_k}^s(1-a_{i_k}^2)\prod\limits_{\substack{1\le \ell\le s\\ \ell\ne k}}e(a_{i_k};a_{i_\ell})}.
%\label{eq:key_c}
\end{gather*}
\end{theorem}

\begin{proof}
Without loss of generality, it suf\/f\/ices to show that
\begin{gather}
C_0\langle e_1e_2\cdots e_s,\xi\rangle
+\sum_{i=1}^s C_i \,\langle e_1\cdots \widehat{e}_i\cdots e_s,\xi\rangle=0,
\label{eq:key0}
\end{gather}
where the coef\/f\/icients $C_0$ and $C_i$ are given by
\begin{gather}
C_0=1-a_1a_2\cdots a_{2s+2}
\qquad\mbox{and}\qquad
C_i=\frac{\prod\limits_{m=1}^{2s+2}(1-a_ia_m)}{a_i^s(1-a_i^2)\prod\limits_{\substack{1\le k\le s\\  k\ne i}}e(a_i;a_k)}.
\label{eq:key0_c}
\end{gather}
Set
$
F(z)=\prod\limits_{m=1}^{2s+2}(a_m-z)$ and $G(z)=\prod\limits_{m=1}^{2s+2}(1-a_mz)
$.
Then, from Lemma \ref{lem:cohom=0}, it follows that
\begin{gather}
\int_0^{\xi\infty}
\!\!\!
\Phi(z)\nabla\left(\frac{F(z)}{z^{s+1}}\right)
\frac{d_q z}{z}=0,\qquad\mbox{where}\qquad \nabla\left(\frac{F(z)}{z^{s+1}}\right)=\frac{F(z)-G(z)}{z^{s+1}}.
\label{eq:nabla(F/z)}
\end{gather}
Since $(F(z)-G(z))/z^{s+1}$ is skew-symmetric under the ref\/lection $z\to z^{-1}$, it is divisible by $z-z^{-1}$,
and we can expand it as
\begin{gather}
\frac{F(z)-G(z)}{z^{s+1}(z-z^{-1})}=
C_0\, e(z;a_1)e(z;a_2)\cdots e(z;a_s)
+\sum_{i=1}^s \! C_i \, e(z;a_1)\cdots \widehat{e}{\hskip 1pt}(z;a_i)\cdots e(z;a_s),\!\!
\label{eq:expand(F-G/z)}
\end{gather}
where the coef\/f\/icients $C_i$ will be determined below.
We obtain $C_0=1-a_1a_2\cdots a_{2s+2}$ from the principal term of asymptotic behavior of (\ref{eq:expand(F-G/z)})
as $z\to +\infty$.
If we put $z=a_i$ $(1\le i\le s)$, then we have
\[
-\frac{F(a_i)-G(a_i)}{a_i^{s}(1-a_i^2)}=C_i  \prod\limits_{\substack{1\le k\le s\\  k\ne i}}e(a_i;a_k).
\]
Since $F(a_i)=0$ and $G(a_i)=\prod\limits_{m=1}^{2s+2}(1-a_ia_m)$ by def\/inition, the above equation implies (\ref{eq:key0_c}).
From (\ref{eq:nabla(F/z)}) and (\ref{eq:expand(F-G/z)}), we obtain (\ref{eq:key0}), which
completes the proof.
\end{proof}

\section[The case $a_1a_2 ... a_{2s+2}\ne 1$]{The case $\boldsymbol{a_1a_2\cdots a_{2s+2}\ne 1}$}

\subsection[$q$-difference equation]{$\boldsymbol{q}$-dif\/ference equation}

Set
\begin{gather}
v_k(z):=\left\{
\begin{array}{ll}
e_{i_1}e_{i_2}\cdots e_{i_{s-1}}(z)&\mbox{if}\quad k=0,\\
e_{i_1}\cdots \widehat{e}_{i_k}\cdots e_{i_{s-1}}(z)
& \mbox{if}\quad    1\le k\le s-1,
\end{array}\right.
\label{eq:uk}
\end{gather}
where the hat symbol denotes the term to be omitted.

Let $T_{a_j}$ be the dif\/ference operator corresponding to the $q$-shift $a_j\to qa_j$.

\begin{theorem}
\label{thm:difference}
%%%%%%%%%%%%%%%%%%%%%%%%%%%%%%
Suppose $a_1a_2{\cdots}a_{2s+2}\!\ne\! 1$.
For the $BC_1$-type Jackson integrals,
if $\{i_1,i_2,\ldots, i_{s-1}\}\!$ $\subset \{1,2,\ldots, 2s+2\}$ and
$j\not\in\{i_1,i_2,\ldots, i_{s-1}\}$, then
the first order vector-valued $q$-difference equation
with respect to the basis $\{v_0,v_1,\ldots,v_{s-1} \}$ defined by {\rm (\ref{eq:uk})}
is given by
\begin{gather}
T_{a_j}(\langle
v_0,\xi
\rangle
,\ldots,
\langle
v_{s-1},\xi
\rangle
)=(\langle
v_0,\xi
\rangle
,\ldots,
\langle
v_{s-1},\xi
\rangle
)B,
\label{eq:difference}
\end{gather}
where $B=UL$. Here $U$ and $L$ are the $s\times s$ matrices  defined by
\[
U=\left(
\begin{array}{ccccc}
c_0&   1&1&\cdots&1   \\
&   c_{1}&&&   \\
&   &c_{2}& &    \\
&   & &\ddots&    \\
&   & & &c_{s-1}
\end{array}
\right) , \qquad
L=
\left(
\begin{array}{lcccc}
 \ 1&   &&&   \\
d_1&
   1&&&   \\
d_2&   &  1 & &    \\
\ \vdots&   & &  \ddots&   \\
d_{s-1}\!\!
&   & & & 1
\end{array}
\right),
\]
where
\[
c_0=\frac{\prod\limits_{m=1}^{2s+2}(1-a_{j}a_m)}{(-a_{j})^s(1-a_1a_2\cdots a_{2s+2})\big(1-a_{j}^2\big)
\prod\limits_{\ell=1}^{s-1}e(a_{i_\ell};a_{j})}
\]
and
\begin{gather}
c_k=e(a_{i_k};a_j)
,\qquad
d_k=\frac{\prod\limits_{m=1}^{2s+2}(1-a_{i_k}a_m)}{(-a_{i_k})^s(1\!-\!a_1a_2\cdots a_{2s+2})(1\!-\!a_{i_k}^2)
e(a_j;a_{i_k})\!\prod\limits_{\substack{1\le \ell\le s-1\\ \ell\ne k}}\!\!e(a_{i_\ell};a_{i_k})}\!\!
\label{eq:cd}
\end{gather}
for $k=1,2,\ldots,s-1$.
Moreover,
\[
\det B=\frac{\prod\limits_{m=1}^{2s+2}(1-a_{j}a_m)}{(-a_{j})^s(1-a_1a_2\cdots a_{2s+2})\big(1-a_{j}^2\big)}.
\]
\end{theorem}

\begin{proof} Equation~(\ref{eq:difference}) is rewritten as
$T_{a_j}(\langle
v_0,\xi
\rangle
,\ldots,
\langle
v_{s-1},\xi
\rangle
)L^{-1}=(\langle
v_0,\xi
\rangle
,\ldots,
\langle
v_{s-1},\xi
\rangle
)U$, where
\[
L^{-1}=\left(\!\!
\begin{array}{lcccc}
\ \ 1&   &&  \\
-d_1&\!\!\!   1&&   \\
-d_2&   &\  1&     \\
\ \ \vdots&   & &\!\ddots\!\!   \\
-d_{s-1}\!\!&   & & &\ 1
\end{array}
\right).
\]
Since $
T_{a_j}\langle
v_i,\xi
\rangle=
\langle
e_jv_i,\xi
\rangle
$,
the above equation is equivalent to
\begin{gather}
\langle
e_jv_0,\xi
\rangle
-\sum_{k=1}^{s-1} d_k
\langle
e_jv_k,\xi
\rangle
=c_0\langle
v_0,\xi
\rangle
\label{eq:Uentry0}
\end{gather}
and
\begin{gather}
\langle
e_jv_k,\xi
\rangle
=\langle
v_0,\xi
\rangle
+c_{k}
\langle
v_k,\xi
\rangle
\qquad\mbox{for}\qquad k=1,2,\ldots, s-1,
\label{eq:Uentryi}
\end{gather}
which are to be proved.
Equation~(\ref{eq:Uentry0}) is a direct consequence of (\ref{eq:e-2}) and Theorem~\ref{thm:key}
if $a_1a_2\cdots a_{2s+2}\ne 1$.
Equation~(\ref{eq:Uentryi}) is trivial using
$
e(z;a_j)=e(z;a_{i_k})+e(a_{i_k};a_j)
$
from (\ref{eq:e-1}). Lastly
$\det B=\det U\det L=c_0 c_1\cdots c_{s-1}$,
which completes the proof.
\end{proof}

Since the function $\Theta(z)$ satisf\/ies $T_{a_j}\Theta(z)=-a_j\Theta(z)$, we immediately have the following
from Theorem \ref{thm:difference}:
%%%%%%%%%%%%%%%%%%%%%%%%%%%%%%
\begin{corollary}
%%%%%%%%%%%%%%%%%%%%%%%%%%%%%%
Suppose $a_1a_2\cdots a_{2s+2}\ne 1$.
For the regularized $BC_1$-type Jackson integrals,
if $\{i_1,i_2,\ldots, i_{s-1}\}\subset \{1,2,\ldots, 2s+2\}$ and
$j\not\in\{i_1,i_2,\ldots, i_{s-1}\}$, then
the first order vector-valued $q$-difference equation
with respect to the basis $\{v_0,v_1,\ldots,v_{s-1}\}$
is given by
\begin{gather}
T_{a_j}(\langle\!\langle
v_0,\xi
\rangle\!\rangle
,\ldots,
\langle\!\langle
v_{s-1},\xi
\rangle\!\rangle
)=(\langle\!\langle
v_0,\xi
\rangle\!\rangle
,\ldots,
\langle\!\langle
v_{s-1},\xi
\rangle\!\rangle
)\bar{B},
\label{eq:eqBbar}
\end{gather}
where
\[
\bar{B}=\frac{-1}{a_j}\left(
\begin{array}{ccccc}
c_0&   1&1&\cdots&1   \\
&   c_{1}&&&   \\
&   &c_{2}& &    \\
&   & &\ddots&    \\
&   & & &c_{s-1}\!\!
\end{array}
\right)\!\!
\left(\!
\begin{array}{lcccc}
\ 1&   &&&   \\
d_1&\!\!\!   1&&&   \\
d_2&   &\  1& &    \\
\ \vdots&   & &\ \ddots\!\!&   \\
d_{s-1}\!\!&   & & &\ 1
\end{array}
\right)
\]
and $c_i$ and $d_i$ are given by {\rm (\ref{eq:cd})}.
In particular, the diagonal entries of the upper triangular part are written as
\[
-\frac{c_0}{a_j}=\frac{\prod\limits_{m=1}^{2s+2}\big(1-a_j^{-1}a_m^{-1}\big)}
{\big(1-a_1^{-1}a_2^{-1}\cdots a_{2s+2}^{-1}\big)\big(1-a_j^{-2}\big)
\prod\limits_{\ell=1}^{s-1}\big(1-a_{i_\ell}a_j^{-1}\big)\big(1-a_{i_\ell}^{-1}a_j^{-1}\big)}
\]
and
\[
-\frac{c_k}{a_j}=\big(1-a_{i_k}a_j^{-1}\big)\big(1-a_{i_k}^{-1}a_j^{-1}\big)
\qquad\mbox{for}\qquad i=1,2,\ldots,s-1.
\]
Moreover,
\begin{gather}
\det \bar{B}=\frac{\prod\limits_{m=1}^{2s+2}\big(1-a_j^{-1}a_m^{-1}\big)}{\big(1-a_j^{-2}\big)\big(1-a_1^{-1}a_2^{-1}\cdots a_{2s+2}^{-1}\big)}.
\label{eq:detBbar}
\end{gather}
\end{corollary}

\begin{remark}
The $q$-dif\/ference system for the $BC_n$-type Jackson integral
is discussed in \cite{AI2} for its rank, and
in \cite{AI3,AI4}
for the explicit expression of the determinant of the coef\/f\/icient matrix of the system.
On the other hand, though it is only for the $BC_1$-type Jackson integral,
the coef\/f\/icient matrix in its Gauss decomposition form is obtained explicitly
only in the present paper.
\end{remark}

\subsection{Application}
The aim of this subsection is to give a simple proof of Gustafson's multiple
$C_n$-type summation formula (Corollary~\ref{cor:Gus}).
The point of the proof is to obtain a recurrence relation of
Gustafson's multiple series of $C_n$-type.
Before we state the recurrence relation, we f\/irst give the def\/inition of the
multiple series of $C_n$-type $\langle\!\langle 1, x \rangle\!\rangle_{\mbox{\tiny\!\rm G}}$.

For $z=(z_1,\ldots,z_n)\in ({\mathbb C}^*)^n$, we set
\begin{gather*}
\Phi_{\mbox{\tiny\!\rm G}}(z):=\prod_{i=1}^n
\prod_{m=1}^{2s+2} z_i^{1/2-\alpha_m}
\frac{(qa_m^{-1}z_i)_\infty}{(a_m z_i)_\infty},
\\
\Delta_{C_n}(z):=\prod_{i=1}^n
\frac{1-z_i^2}{z_i}\prod_{1\le j<k \le n}
\frac{(1-z_j/z_k)(1-z_jz_k)}{z_j},
\end{gather*}
where
$
q^{\alpha_m}=a_m
$.
For an arbitrary $\xi=(\xi_1,\ldots,\xi_n)\in ({\mathbb C}^*)^n$,
we def\/ine the $q$-shift $\xi\to q^\nu\xi$ by a lattice point $\nu=(\nu_1,\ldots,\nu_n)\in {\mathbb Z}^n$,
where
$q^\nu \xi:=(q^{\nu_1}\xi_1,\ldots,q^{\nu_n}\xi_n)\in ({\mathbb C}^*)^n$.
For $\xi=(\xi_1,\ldots,\xi_n)\in ({\mathbb C}^*)^n$
we def\/ine the sum over the lattice ${\mathbb Z}^n$ by
\begin{gather*}
\langle 1,\xi\rangle_{\mbox{\tiny\!\rm G}}:=
(1-q)^n\sum_{\nu\in {\mathbb Z}^n}
\Phi_{\mbox{\tiny\!\rm G}}(q^\nu \xi)\Delta_{C_n}(q^\nu \xi),
\end{gather*}
which we call the {\it $BC_n$-type Jackson integral}.
Moreover we set $\langle\!\langle 1, \xi \rangle\!\rangle_{\mbox{\tiny\!\rm G}}:=\langle 1, \xi \rangle_{\mbox{\tiny\!\rm G}}/\Theta_{\mbox{\tiny\!\rm G}}(\xi)$,
where
\begin{gather*}
\Theta_{\mbox{\tiny\!\rm G}}(\xi)
:=\prod_{i=1}^n
\frac{\xi_i^{i-\alpha_1-\alpha_2-\cdots-\alpha_{2s+2}}\theta(\xi_i^2)}
         {\prod\limits_{m=1}^{2s+2}\theta(a_m \xi_i)}
\prod_{1\le j<k \le n}\!\!\!\theta(\xi_j/\xi_k)\theta(\xi_j\xi_k).
%\label{eq:Theta}
\end{gather*}
By def\/inition, it can be conf\/irmed that $\langle\!\langle 1, \xi \rangle\!\rangle_{\mbox{\tiny\!\rm G}}$ is holomorphic on $({\mathbb C}^*)^n$
(see \cite[Proposition~3.7]{AI3}), and
we call it the {\it regularized $BC_n$-type Jackson integral}.
In particular, if we assume $s=n$ we call $\langle\!\langle 1, \xi \rangle\!\rangle_{\mbox{\tiny\!\rm G}}$ the
{\it regularized Jackson integral of Gustafson's $C_n$-type},
which is, in particular, a~constant not depending on $\xi\in ({\mathbb C^*})^n$.

\begin{remark}
For further results on $BC_n$-type Jackson integrals,
see \cite{AI1,AI3,AI2,AI4,Ito9,IO}, for instance.
\end{remark}

Now we state the recurrence relation for Gustafson's sum $\langle\!\langle 1, \xi \rangle\!\rangle_{\mbox{\tiny\!\rm G}}$.
%%%%%%%%%%%%%%%%%%%%%%%%%%%%%%%%%%
\begin{proposition}
\label{prop:recGus1}
Suppose $s=n$ and $x = (x_1, \ldots, x_n) \in ({\mathbb C}^*)^n$. The sum $\langle\!\langle 1, x \rangle\!\rangle_{\mbox{\tiny\!\rm G}}$ satisfies
\[
T_{a_j}\langle\!\langle 1, x \rangle\!\rangle_{\mbox{\tiny\!\rm G}}=\langle\!\langle 1, x \rangle\!\rangle_{\mbox{\tiny\!\rm G}} \frac{\prod\limits_{m=1}^{2n+2}(1-a_j^{-1}a_m^{-1})}{\big(1-a_j^{-2}\big)\big(1-a_1^{-1}a_2^{-1}\cdots a_{2n+2}^{-1}\big)}
\qquad\mbox{for}\qquad j=1,2,\ldots,2n+2.
\]
\end{proposition}
%%%%%%%%%%%%%%%%%%%%%%%%%%%%%%
%
\begin{proof} We assume $s=n$ for the basis $\{v_0,v_1,\ldots,v_{n-1}\}$ of the $BC_1$-type Jackson integral.
Let~$P$ be the transition matrix from the basis $\{v_0,v_1,\ldots,v_{n-1}\}$ to
$\{\chi_{(n-1)},\chi_{(n-2)},\ldots,\chi_{(0)}\}$:
\[
(\chi_{(n-1)},\chi_{(n-2)},\ldots,\chi_{(0)})=(v_0,v_1,\ldots,v_{n-1})P,
\]
where $\chi_{(i)}$ is the irreducible character of type $C_1$ def\/ined by
\begin{gather}
\chi_{(i)}(z)=\frac{z^{i+1}-z^{-i-1}}{z-z^{-1}} \qquad\mbox{for}\qquad i=0,1,2,\ldots.
\label{eq:chi}
\end{gather}
From (\ref{eq:eqBbar}) it follows that
\begin{gather}
T_{a_j}(\langle\!\langle
\chi_{(n-1)},\xi
\rangle\!\rangle
,\ldots,
\langle\!\langle
\chi_{(0)}
,\xi
\rangle\!\rangle
)=(\langle\!\langle
\chi_{(n-1)},\xi
\rangle\!\rangle
,\ldots,
\langle\!\langle
\chi_{(0)},\xi
\rangle\!\rangle
)P^{-1}\bar{B}P,
\label{eq:PBbar/P}
\end{gather}
so that
\begin{gather}
T_{a_j}\det \big( \langle\!\langle \chi_{(n-i)}, x_j \rangle\!\rangle \big)_{\!1 \le i, j \le n}
=
\det \big( \langle\!\langle \chi_{(n-i)}, x_j \rangle\!\rangle \big)_{\!1 \le i, j \le n}\det \bar{B}
\label{eq:TdetBbar}
\end{gather}
for $x=(x_1,\ldots,x_n)\in {({\mathbb C}^*)}^n$.
By def\/inition, the relation between
the determinant of the $BC_1$-type Jackson integrals and the Jackson integral of Gustafson's $C_n$-type itself
is given as
\begin{gather}
\det \big( \langle\!\langle \chi_{(n-i)}, x_j \rangle\!\rangle \big)_{\!1 \le i, j \le n}
 =
\langle\!\langle 1,x \rangle\!\rangle_{\mbox{\tiny\!\rm G}}
\prod_{1 \le j < k \le n}\frac{ \theta(x_j/x_k) \theta(x_j x_k) }{ x_j },
\label{eq:<<BC1-BCn>>}
\end{gather}
which is also referred to in \cite{IO}.
From (\ref{eq:TdetBbar}) and (\ref{eq:<<BC1-BCn>>}), we obtain $T_{a_j}\langle\!\langle 1, x \rangle\!\rangle_{\mbox{\tiny\!\rm G}}
=\langle\!\langle 1, x \rangle\!\rangle_{\mbox{\tiny\!\rm G}} \det \bar{B}$, where $\det\bar{B}$ has already been given in (\ref{eq:detBbar}).
\end{proof}
\begin{remark}
 The explicit form of the coef\/f\/icient matrix of the system~(\ref{eq:PBbar/P}) is given in~\cite{A1} or~\cite{AI2}.
\end{remark}
%%%%%%%%%%%%%%%%%%%%%%%%%%%%%%%%%
\begin{corollary}[Gustafson \cite{Gu4}]
\label{cor:Gus}
%%%%%%%%%%%%%%%%%%%%%%%%%%%%%%%%%
Suppose $s=n$ and $x = (x_1, \ldots, x_n) \in ({\mathbb C}^*)^n$.
Then the sum $\langle\!\langle 1, x \rangle\!\rangle_{\mbox{\tiny\!\rm G}}$ is written as
\[
\langle\!\langle 1, x \rangle\!\rangle_{\mbox{\tiny\!\rm G}}
 =
(1-q)^n
\frac{ (q)_\infty^n\prod\limits_{1 \le i < j \le 2n+2} \big(q a_i^{-1} a_j^{-1}\big)_\infty }
     { \big(q a_1^{-1} a_2^{-1} \cdots a_{2n+2}^{-1}\big)_\infty }.
\]
\end{corollary}

\begin{proof}
By repeated use of the recurrence relation in Proposition \ref{prop:recGus1},
using the asymptotic behavior of the Jackson integral as the boundary condition of the recurrence relation,
we eventually obtain Corollary~\ref{cor:Gus}.
See \cite{Ito8} for further details about the proof. \end{proof}
\begin{remark}
From (\ref{eq:<<BC1-BCn>>}) and Corollary \ref{cor:Gus}, we see
\begin{gather*}
\det \big( \langle\!\langle \chi_{(s-i)}, x_j \rangle\!\rangle \big)_{\!1 \le i, j \le s}
=(1-q)^s
\frac{ (q)_\infty^s\prod\limits_{1 \le i < j \le 2s+2} \big(q a_i^{-1} a_j^{-1}\big)_\infty }
{ \big(q a_1^{-1} a_2^{-1} \cdots a_{2s+2}^{-1}\big)_\infty }
\prod_{1 \le j < k \le s}\frac{ \theta(x_j/x_k) \theta(x_j x_k) }{ x_j },
\end{gather*}
which is non-degenerate under generic condition. This indicates that
the set $\{\chi_{(s-1)},\chi_{(s-2)},\ldots$, $\chi_{(0)}\}$ is linearly independent.
And we eventually know the rank of the $q$-dif\/ference system with respect to this basis is $s$,
so are the ranks of the systems (\ref{eq:difference}) and (\ref{eq:eqBbar}).
\end{remark}

\section[The case $a_1a_2 ... a_{2s+2}= 1$]{The case $\boldsymbol{a_1a_2\cdots a_{2s+2}= 1}$}

\subsection[Reflection equation]{Ref\/lection equation}

\begin{theorem}
\label{thm:ref-eqn}
%%%%%%%%%%%%%%%%%%%%%%%%%%%%%%
Suppose $a_1a_2\cdots a_{2s+2}= 1$.
Let $v_k(z)$, $k=1,2,\ldots, s-1$, be the functions defined by {\rm (\ref{eq:uk})} for the fixed indices $i_1,i_2,\ldots,i_{s-1}\in \{1,2,\ldots,2s+2\}.$
If $j_1,j_2\not\in\{i_1,i_2,\ldots,i_{s-1}\}$, then
\[
(\langle e_{j_1}v_1,\xi\rangle,\ldots,\langle e_{j_1}v_{s-1},\xi\rangle)=
(\langle e_{j_2}v_1,\xi\rangle,\ldots,\langle e_{j_2}v_{s-1},\xi\rangle)M,
\]
where $M=M_{j_2}NM_{j_1}^{-1}$. Here $M_j$ and $N$ are the matrices defined by
\[
M_j=
\left(\!\!
\begin{array}{lccccc}
\gamma_{1,j}&   &&&   \\
\gamma_{2,j}&\!\!\!   1&&&   \\
\gamma_{3,j}&   &\  1& &    \\
\ \vdots&   & &\!\ddots\!\!&   \\
\gamma_{s-1,j}\!\!&   & & &\ \ 1
\end{array}
\right),\qquad
N=
\left(\!\!
\begin{array}{cccccc}
1&  \sigma_{2}&\sigma_{3}&\cdots&\sigma_{s-1}   \\
&   \tau_{2}&&&   \\
&   &\tau_{3}& &    \\
&   & &\ddots&    \\
&   & & &\tau_{s-1}\!\!\!
\end{array}
\right),
\]
where the entries of the above matrices are given by
\begin{gather*}
\sigma_k=\frac{e(a_{j_1};a_{j_2})}{e(a_{i_k};a_{j_2})},\qquad
\tau_k=\frac{e(a_{j_1};a_{i_k})}{e(a_{j_2};a_{i_k})},\\
\gamma_{k,j}=\frac{a_j^s\big(1-a_j^2\big)}{a_{i_k}^s\big(1-a_{i_k}^2\big)}\prod_{m=1}^{2s+2}\frac{1-a_{i_k}a_m}{1-a_{j}a_m}
\prod\limits_{\substack{1\le \ell\le s-1\\  \ell\ne k}}\frac{e(a_j;a_{i_\ell})}{e(a_{i_k};a_{i_\ell})}.
\end{gather*}
Moreover,
\begin{gather}
\det M=
\frac{a_{j_2}^s\big(1-a_{j_2}^2\big)}{a_{j_1}^s\big(1-a_{j_1}^2\big)}
\prod_{m=1}^{2s+2}\frac{1-a_{j_1}a_m}{1-a_{j_2}a_m}.
\label{eq:detM}
\end{gather}
\end{theorem}
%%%%%%%%%%%%%%%%%%%%%%%%%%%%%%
%
\begin{proof}
First we will prove the following:
\begin{gather}
(\langle e_{j_1}v_1,\xi\rangle,\langle e_{j_1}v_2,\xi\rangle,\ldots,\langle e_{j_1}v_{s-1},\xi\rangle)M_{j_1}=
(\langle v_0,\xi\rangle,\langle e_{j_2}v_2,\xi\rangle,\ldots,\langle e_{j_2}v_{s-1},\xi\rangle)N,
\label{eq:MN1}
\\
(\langle e_{j_2}v_1,\xi\rangle,\langle e_{j_2}v_2,\xi\rangle,\ldots,\langle e_{j_2}v_{s-1},\xi\rangle)M_{j_2}=
(\langle v_0,\xi\rangle,\langle e_{j_2}v_2,\xi\rangle,\ldots,\langle e_{j_2}v_{s-1},\xi\rangle),
\label{eq:MN2}
\end{gather}
which are equivalent to
\begin{gather}
\sum_{k=1}^{s-1}\gamma_{k,j} \langle e_j v_k,\xi\rangle=\langle v_0,\xi\rangle
\label{eq:Mentry0}
\end{gather}
and
\begin{gather}
\langle e_{j_1}v_k,\xi\rangle=\sigma_k \langle v_0,\xi\rangle+\tau_k \langle e_{j_2}v_k,\xi\rangle,
\qquad k=2,\ldots, s-1.
\label{eq:Mentryi}
\end{gather}
Under the condition $a_1a_2\cdots a_{2s+2}= 1$,
Equation~(\ref{eq:Mentry0}) is a direct consequence of Theorem \ref{thm:key}.
Equation~(\ref{eq:Mentryi}) is trivial from the equation
\[
e(z;a_i)=e(z;a_j)\frac{e(a_i;a_k)}{e(a_j;a_k)}+e(z;a_k)\frac{e(a_i;a_j)}{e(a_k;a_j)},
\]
which was given in (\ref{eq:e-3}).
From (\ref{eq:MN1}) and (\ref{eq:MN2}), it follows
$
M=M_{j_2}NM_{j_1}^{-1}
$.
Moreover, we obtain
\[\det M=\frac{\det M_{j_1}\det N}{\det M_{j_2}}=
\frac{\gamma_{1,j_1}\tau_2\cdots \tau_{s-1}}{\gamma_{1,j_2}}=
\frac{a_{j_2}^s\big(1-a_{j_2}^2\big)}{a_{j_1}^s\big(1-a_{j_1}^2\big)}
\prod_{m=1}^{2s+2}\frac{1-a_{j_1}a_m}{1-a_{j_2}a_m},\]
which completes the proof.\end{proof}
%%%%%%%%%%%%%%%%%%%%%%%%%%%%%%
\begin{corollary}
\label{cor:s=2Taj}
%%%%%%%%%%%%%%%%%%%%%%%%%%%%%%
Suppose $s=2$ and  the condition $a_6=\frac{q}{a_1a_2a_3a_4a_5}$.
The recurrence relation for the $BC_1$-type Jackson integral $\langle 1,\xi\rangle$ is
\[
T_{a_j}\langle 1,\xi\rangle=\langle 1,\xi\rangle\frac{q}{a_ja_6}
\prod_{\substack{1\le \ell\le 5 \\ \ell\ne j}}\frac{1-a_ja_\ell}{1-qa_6^{-1}a_\ell^{-1}}
\qquad\mbox{for}\qquad j=1,2,\ldots, 5.
\]
\end{corollary}
%%%%%%%%%%%%%%%%%%%%%%%%%%%%%%
%
\begin{proof}
Without loss of generality, it suf\/f\/ices to show that
\begin{gather}
T_{a_1}\langle 1,\xi\rangle=\langle 1,\xi\rangle\frac{q}{a_1a_6}\prod_{\ell=2}^5\frac{1-a_1a_\ell}{1-qa_6^{-1}a_\ell^{-1}}.
\label{eq:Ta1<1,xi>}
\end{gather}
Set
$J(a_1,a_2,a_3,a_4,a_5,a_6;\xi):=\langle 1,\xi\rangle$.
Under the condition $a_1a_2a_3a_4a_5a_6=1$, we have
\[
J(qa_1,a_2,a_3,a_4,a_5,a_6;\xi)=J(a_1,a_2,a_3,a_4,a_5,qa_6;\xi)
\frac{a_6^2}{a_1^2}\prod_{\ell=2}^5\frac{1-a_1a_\ell}{1-a_6a_\ell}
\]
from Theorem \ref{thm:ref-eqn} by setting $j_1=1$ and $j_2=6$.
We now replace $a_6$ by $q^{-1}a_6$ in the above equation.
Then, under the condition $a_1a_2a_3a_4a_5(q^{-1}a_6)=1$,
we have
\[
J\left(qa_1,a_2,a_3,a_4,a_5,q^{-1}a_6;\xi\right)=J(a_1,a_2,a_3,a_4,a_5,a_6;\xi)
\frac{q}{a_1a_6}\prod_{\ell=2}^5\frac{1-a_1a_\ell}{1-qa_6^{-1}a_\ell^{-1}}.
\]
Since
$T_{a_1}\langle 1,\xi\rangle=J(qa_1,a_2,a_3,a_4,a_5,q^{-1}a_6;\xi)$
under this condition $a_6=q(a_1a_2a_3a_4a_5)^{-1}$, we obtain (\ref{eq:Ta1<1,xi>}), which completes the proof. \end{proof}
%%%%%%%%%%%%%%%%%%%%%%%%%%%%%%
\begin{corollary}
\label{cor:s=n+1Taj}
%%%%%%%%%%%%%%%%%%%%%%%%%%%%%%
Suppose $s=n+1$ and  the condition $a_{2n+4}=\frac{q}{a_1a_2\cdots a_{2n+3}}$. Then
the recurrence relation for Gustafson's sum $\langle 1,x\rangle_{\mbox{\tiny\!\rm G}}$ where $x=(x_1,\ldots,x_n)\in ({\mathbb C}^*)^n$ is given by
\[
T_{a_j}\langle 1,x\rangle_{\mbox{\tiny\!\rm G}}=\langle 1,x\rangle_{\mbox{\tiny\!\rm G}}\frac{q}{a_ja_{2n+4}}
\prod_{\substack{1\le \ell\le 2n+3 \\ \ell\ne j}}\frac{1-a_ja_\ell}{1-qa_\ell^{-1}a_{2n+4}^{-1}}
\qquad\mbox{for}\qquad j=1,2,\ldots, 2n+3.
\]
\end{corollary}
%%%%%%%%%%%%%%%%%%%%%%%%%%%%%%
%
\begin{proof} Fix $s=n+1$. For the $BC_1$-type Jackson integral, we f\/irst set
\[
J(a_1,a_2,\ldots,a_{2n+4};x):=\det \big( \langle \chi_{(n-i)}, x_j \rangle \big)_{\!1 \le i, j \le n},
\]
where $\chi_{(i)}$ is def\/ined in (\ref{eq:chi}),
under no condition on $a_1,a_2,\ldots,a_{2n+4}$. By the def\/inition of $\Phi$, we have
\[
J(qa_1,a_2,\ldots,a_{2n+4};x)=\det \big( \langle e_1\chi_{(n-i)}, x_j \rangle \big)_{\!1 \le i, j \le n}.
\]
Let $Q$ be the transition matrix from the basis
$\{v_1,v_2,\ldots,v_n\}$ to $\{\chi_{(n-1)},\chi_{(n-2)},\ldots,\chi_{(0)}\}$, i.e.,
\[
(\chi_{(n-1)},\chi_{(n-2)},\ldots,\chi_{(0)})=(v_1,v_2,\ldots,v_n)Q.
\]
Under the condition  $a_1a_2\cdots a_{2n+4}=1$,
from Theorem \ref{thm:ref-eqn} with $j_1=1$ and $j_2=2n+4$, it follows that
\[
(\langle e_1v_1,\xi\rangle,\ldots,\langle e_1v_{n},\xi\rangle)=
(\langle e_{2n+4}v_1,\xi\rangle,\ldots,\langle e_{2n+4}v_{n},\xi\rangle)M,
\]
so that
\[
(\langle e_1\chi_{(n-1)},\xi\rangle,\ldots,\langle e_1\chi_{(0)},\xi\rangle)=
(\langle e_{2n+4}\chi_{(n-1)},\xi\rangle,\ldots,\langle e_{2n+4}\chi_{(0)},\xi\rangle)Q^{-1}MQ.
\]
This indicates that
\[
\det \big( \langle e_1\chi_{(n-i)}, x_j \rangle \big)_{\!1 \le i, j \le n}=
\det \big( \langle e_{2n+4}\chi_{(n-i)}, x_j \rangle \big)_{\!1 \le i, j \le n}\det M.
\]
From (\ref{eq:detM}) and the above equation we have
\[J(qa_1,a_2,\ldots,a_{2n+4};x)=J(a_1,a_2,\ldots,qa_{2n+4};x)
\left(\frac{a_{2n+4}}{a_1}\right)^{\!\!n+1}\prod_{\ell=2}^{2n+3}\frac{1-a_\ell a_1}{1-a_\ell a_{2n+4}},\]
under the condition  $a_1a_2\cdots a_{2n+4}=1$.
We now replace $a_{2n+4}$ by $q^{-1}a_{2n+4}$ in the above equation.
Then we have
\[J\left(qa_1,a_2,\ldots,q^{-1}a_{2n+4};x\right)=J(a_1,a_2,\ldots,a_{2n+4};x)
\frac{q}{a_1a_{2n+4}}\prod_{\ell=2}^{2n+3}\frac{1-a_\ell a_1}{1-qa_\ell^{-1} a_{2n+4}^{-1}},\]
under the condition $a_1a_2\cdots(q^{-1} a_{2n+4})=1$.
Since \[T_{a_1}\det \big( \langle \chi_{(n-i)}, x_j \rangle \big)_{\!1 \le i, j \le n}
=J\left(qa_1,a_2,\ldots,q^{-1}a_{2n+4};x\right)\]
if $a_1a_2\cdots a_{2n+4}=q$, we have
\[
T_{a_1}\det \big( \langle \chi_{(n-i)}, x_j \rangle \big)_{\!1 \le i, j \le n}
=\det \big( \langle \chi_{(n-i)}, x_j \rangle \big)_{\!1 \le i, j \le n}\,
\frac{q}{a_1a_{2n+4}}\prod_{\ell=2}^{2n+3}\frac{1-a_\ell a_1}{1-qa_\ell^{-1} a_{2n+4}^{-1}}.
\]
On the other hand, if $x=(x_1,\ldots,x_n)\in ({\mathbb C}^*)^n$, then, by def\/inition we have
\begin{gather*}
\det \big( \langle \chi_{(n-i)}, x_j \rangle \big)_{\!1 \le i, j \le n}
 =
\langle 1,x\rangle_{\mbox{\tiny\!\rm G}},
%\label{eq:<BC1-BCn>}
\end{gather*}
which is also referred  to in \cite{IO}.
Therefore, under the condition $a_{2n+4}=q(a_1a_2\cdots a_{2n+3})^{-1}$ we obtain
\[
T_{a_1}\langle 1,x\rangle_{\mbox{\tiny\!\rm G}}
=\langle 1,x\rangle_{\mbox{\tiny\!\rm G}}
\frac{q}{a_1a_{2n+4}}\prod_{\ell=2}^{2n+3}\frac{1-a_\ell a_1}{1-qa_\ell^{-1} a_{2n+4}^{-1}}.
\]
Since the same argument holds for parameters $a_2,\ldots,a_{2n+3}$,
we can conclude Corollary \ref{cor:s=n+1Taj}. \end{proof}
\begin{remark}
If we take $\xi=a_i$, $i=1,\ldots, 6$, and add the terminating condition $a_1a_2=q^{-N}$, $N=1,2,\ldots,$ to the assumptions of Corollary \ref{cor:s=2Taj}, then the f\/inite product expression
of $\langle 1,\xi\rangle$, which is equivalent to Jackson's formula for
terminating ${}_8\phi_7$ series  \cite[p.~43, equation~(2.6.2)]{GR},
is obtained from f\/inite repeated use of Corollary \ref{cor:s=2Taj}. In the same way, if we take a suitable $x$
and add the terminating condition to the
assumptions of Corollary \ref{cor:s=n+1Taj}, then the f\/inite product expression
of $\langle 1,x\rangle_{\mbox{\tiny\!\rm G}}$,
which is equivalent to the Jackson type formula for
terminating multiple~${}_8\phi_7$ series
(see \cite[Theorem 4]{DG} or \cite[p.~231, equation~(4.4)]{DS}, for instance),
is obtained from f\/inite repeated use of Corollary~\ref{cor:s=n+1Taj}.
\end{remark}

\subsection{Application}
The aim of this subsection is to give a simple proof of the following propositions
proved by Nassrallah and Rahman~\cite{NR} and Gustafson~\cite{Gu3}.
%%%%%%%%%%%%%%%%%%%%%%%%%%%%%%
\begin{proposition}[Nassrallah--Rahman]
%%%%%%%%%%%%%%%%%%%%%%%%%%%%%%
Assume $|a_i|<1$ for $1\le i\le 5$. If $a_6=\frac{q}{a_1a_2a_3a_4a_5}$, then
\begin{gather}
\frac{1}{2\pi\sqrt{-1}}\int_{{\mathbb T}}
\frac{\big(qa_6^{-1}z\big)_\infty\big(qa_6^{-1}z^{-1}\big)_\infty\big(z^2\big)_\infty\big(z^{-2}\big)_\infty}
{\prod\limits_{i=1}^5(a_iz)_\infty\big(a_iz^{-1}\big)_\infty}\frac{dz}{z}=
\frac{2\prod\limits_{k=1}^5\big(qa_6^{-1}a_k^{-1}\big)_\infty}{(q)_\infty\prod\limits_{1\le i<j\le 5}(a_ia_j)_\infty},
\label{eq:NR}
\end{gather}
where ${\mathbb T}$ is the unit circle taken in the positive direction.
\end{proposition}
%%%%%%%%%%%%%%%%%%%%%%%%%%%%%%
%
\begin{proof} We denote the left-hand side of (\ref{eq:NR}) by $I(a_1,a_2,a_3,a_4,a_5)$. By residue calculation,
\begin{gather}
I(a_1,a_2,a_3,a_4,a_5)
 =\sum_{k=1}^5\sum_{\nu=0}^\infty\mathop{\hbox{\rm Res }}_{z=a_k q^\nu}\!
\left[\frac{\theta\big(qa_6^{-1}z^{-1}\big)\theta\big(z^{-2}\big)}{z\prod\limits_{m=1}^5\theta\big(a_m z^{-1}\big)}
%\times
z(1-z^2)
\prod_{m=1}^6
\frac{\big(qa_m^{-1}z\big)_\infty}{(a_m z)_\infty}\right]\frac{dz}{z}\!\!
\label{eq:Res1}\\
\phantom{I(a_1,a_2,a_3,a_4,a_5)}{} =\sum_{k=1}^5
\left[\mathop{\hbox{\rm Res }}_{z=a_k}\frac{\theta\big(qa_6^{-1}z^{-1}\big)\theta\big(z^{-2}\big)}{z\prod\limits_{m=1}^5\theta\big(a_m z^{-1}\big)}
\frac{dz}{z}\right]
\int_0^{a_k\infty}\!\!\!
%\hskip -10pt
z(1-z^2)
\prod_{m=1}^6
\frac{\big(qa_m^{-1}z\big)_\infty}{(a_m z)_\infty}
\frac{d_qz}{z}\nonumber\\
\phantom{I(a_1,a_2,a_3,a_4,a_5)}{}  =\sum_{k=1}^5 R_k \langle 1,a_k\rangle,
\label{eq:Res2}
\end{gather}
where
\[
R_k:=\mathop{\hbox{\rm Res }}_{z=a_k}\frac{\theta\big(qa_6^{-1}z^{-1}\big)\theta\big(z^{-2}\big)}{z\prod\limits_{m=1}^5\theta\big(a_m z^{-1}\big)}\frac{dz}{z}
=
\frac{\theta\big(qa_6^{-1}a_k^{-1}\big)\theta\big(a_k^{-2}\big)}{(q)_\infty^2 a_k \prod\limits_{\substack{1\le m\le 5\\ m\ne k}}\theta\big(a_m a_k^{-1}\big)},
\]
whose recurrence relation is
\begin{gather}
T_{a_j}R_k=\left(q^{-1}{a_ja_6}\right)R_k
\label{eq:TRk}
\end{gather}
for $1\le j,k\le 5$, which is obtained using (\ref{eq:thProperty}).
From (\ref{eq:Res2}), (\ref{eq:TRk}) and Corollary \ref{cor:s=2Taj},
we obtain the recurrence relation for $I(a_1,a_2,a_3,a_4,a_5)$ as
\[
T_{a_j}I(a_1,a_2,a_3,a_4,a_5)=I(a_1,a_2,a_3,a_4,a_5)
\prod_{\substack{1\le\ell\le 5\\ \ell\ne j}}\frac{1-a_ja_\ell}{1-qa_6^{-1}a_\ell^{-1}}.
\]
By repeated use of the above relation, we obtain
\begin{gather*}
I(a_1,a_2,a_3,a_4,a_5)
=
\frac{\prod\limits_{k=1}^5\big(qa_6^{-1}a_k^{-1}\big)_{2N}}{\prod\limits_{1\le i<j\le 5}(a_ia_j)_{2N}}
I\big(q^Na_1,q^Na_2,q^Na_3,q^Na_4,q^Na_5\big) \\
\phantom{I(a_1,a_2,a_3,a_4,a_5)}{} =
\frac{\prod\limits_{k=1}^5\big(qa_6^{-1}a_k^{-1}\big)_{\infty}}{\prod\limits_{1\le i<j\le 5}(a_ia_j)_{\infty}}
\lim_{N\to \infty}I\big(q^Na_1,q^Na_2,q^Na_3,q^Na_4,q^Na_5\big)
\end{gather*}
and
\[
\lim_{N\to \infty}I\big(q^Na_1,q^Na_2,q^Na_3,q^Na_4,q^Na_5\big)=
\frac{1}{2\pi\sqrt{-1}}\int_{{\mathbb T}}\big(z^2\big)_\infty\big(z^{-2}\big)_\infty\frac{dz}{z}=\frac{2}{(q)_\infty}.
\]
This completes the proof. \end{proof}

\begin{remark}
 Strictly speaking, the residue calculation (\ref{eq:Res1}) requires that
\begin{gather}
I_\varepsilon:=
{1\over 2\pi \sqrt{-1}}\int_{|z|=\varepsilon}
\!\!\!\!\!\!
\frac{\big(qa_6^{-1}z\big)_\infty\big(qa_6^{-1}z^{-1}\big)_\infty\big(z^2\big)_\infty\big(z^{-2}\big)_\infty}
{\prod\limits_{i=1}^5(a_iz)_\infty\big(a_iz^{-1}\big)_\infty}\frac{dz}{z}\ \to 0 \qquad\mbox{if}\qquad \varepsilon \to 0,
\label{eq:to0}
\end{gather}
which can be shown in the following way.
We f\/irst take $\varepsilon=q^N\varepsilon'$ for $\varepsilon'>0$ and positive inte\-ger~$N$.
If we put
\[
F(z):=\frac{\big(qa_6^{-1}z\big)_\infty\big(qa_6^{-1}z^{-1}\big)_\infty\big(z^2\big)_\infty\big(z^{-2}\big)_\infty}
{\prod\limits_{i=1}^5(a_iz)_\infty\big(a_iz^{-1}\big)_\infty},
\]
then we have
$F(z)=zG_1(z)G_2(z)$,
where
\[
G_1(z)=\frac{\theta\big(qa_6^{-1}z^{-1}\big)\theta\big(z^{-2}\big)}
{z\prod\limits_{i=1}^5\theta\big(a_i z^{-1}\big)},
\qquad
G_2(z)=\big(1-z^2\big)\prod\limits_{i=1}^6\frac{\big(qa_i^{-1}z\big)_\infty}{({a_i}z)_\infty}.
\]
Since $G_1(z)$ is a continuous function on the compact set $|z|=\varepsilon'$
and is invariant under the $q$-shift $z\to qz$ under the condition $a_6=q(a_1a_2a_3a_4a_5)^{-1}$, $|G_1(z)|$ is bounded on
$|z|=q^N\varepsilon'$.
$|G_2(z)|$ is also bounded because $G_2(z)\to 1$ if $z\to 0$.
Thus there exists $C>0$ such that $|F(z)|<C|z|$.
If we put $z=\varepsilon e^{2\pi\sqrt{-1}\tau}$, then
\[
|I_\varepsilon|<\int_0^1|F(\varepsilon e^{2\pi\sqrt{-1}\tau})|d\tau
<C\int_0^1|\varepsilon e^{2\pi\sqrt{-1}\tau}|d\tau
=C\varepsilon
  \to 0,\qquad  \varepsilon\to 0,
\]
which proves (\ref{eq:to0}).
\end{remark}
%
%%%%%%%%%%%%%%%%%%%%%%%%%%%%%%
\begin{proposition}[Gustafson \cite{Gu3}]
%%%%%%%%%%%%%%%%%%%%%%%%%%%%%%
Assume $|a_i|<1$ for $1\le i\le 2n+3$. If $ a_{2n+4}=\frac{q}{a_1a_2\cdots a_{2n+3}}$, then
\begin{gather}
\left(\frac{1}{2\pi\sqrt{-1}}\right)^{\!\!n}\int_{{\mathbb T}^n}
\prod\limits_{i=1}^n
\frac{\big(qa_{2n+4}^{-1}z_i\big)_\infty\big(qa_{2n+4}^{-1}z_i^{-1}\big)_\infty\big(z_i^2\big)_\infty
\big(z_i^{-2}\big)_\infty}
{\prod\limits_{k=1}^{2n+3}(a_kz_i)_\infty\big(a_kz_i^{-1}\big)_\infty}
\nonumber\\
\qquad \quad {}
\times\prod_{1\le i<j \le n}(z_iz_j)_\infty\big(z_iz_j^{-1}\big)_\infty\big(z_i^{-1}z_j\big)_\infty \big(z_i^{-1}z_j^{-1}\big)_\infty
\frac{dz_1}{z_1}\wedge\cdots\wedge\frac{dz_n}{z_n}
\nonumber\\
\qquad {} =
\frac{2^nn!\prod\limits_{k=1}^{2n+3}\big(qa_{2s+4}^{-1}a_k^{-1}\big)_\infty}{(q)_\infty^n\prod\limits_{1\le i<j\le 2s+3}(a_ia_j)_\infty},
\label{eq:NR-G}
\end{gather}
where ${\mathbb T}^n$ is the $n$-fold direct product of the unit circle traversed in the positive direction.
\end{proposition}

The proof below is based on an idea using residue computation
due to Gustafson \cite{Gu4},
which is done for the case of the hypergeometric
integral under no balancing condition.
Here we will show that his residue method is still ef\/fective even
for the integral under the balancing condition
$a_1a_2\cdots a_{2n+4}=q$.
In particular, this is dif\/ferent from his proof in~\cite{Gu3}.

\begin{proof}
Let $L$ be the set of indices def\/ined by
\[L:=\{\lambda=(\lambda_1,\ldots,\lambda_n) ; \; 1\le \lambda_1<\lambda_2<\cdots<\lambda_n\le 2n+3\}.\]
Set $a_{(\mu)}:=(a_{\mu_1},\ldots,a_{\mu_n})\in ({\mathbb C}^*)^n$ for $\mu=(\mu_1,\ldots,\mu_n)\in L$.
We denote the left-hand side of~(\ref{eq:NR-G}) by $I(a_1,a_2,\ldots,a_{2n+3})$. By residue calculation,
we have
\begin{gather}
I(a_1,a_2,\ldots,a_{2n+3})=\sum_{\mu\in L}R_{\mu}\langle 1,a_{(\mu)}\rangle_{\mbox{\tiny\!\rm G}},
\label{eq:Res}
\end{gather}
where the coef\/f\/icients $R_\mu$, $\mu\in L$, are
\[
R_\mu:=
\mathop{\hbox{\rm Res }}_{\substack{z_1=a_{\mu_1}\\ \cdots \\ z_n=a_{\mu_n}}}\left[
\prod_{i=1}^n
\frac{\theta\big(qa_{2n+4}^{-1}z_i^{-1}\big)\theta\big(z_i^{-2}\big)}{z_i\prod\limits_{m=1}^{2n+3}\theta\big(a_m z_i^{-1}\big)}
\prod_{1\le j<k\le n}\theta\big(z_j^{-1}z_k\big)\theta\big(z_j^{-1}z_k^{-1}\big)\right]
\frac{dz_1}{z_1}\wedge\cdots\wedge\frac{dz_n}{z_n}.
\]
The recurrence relation for $R_\mu$ is
\begin{gather}
T_{a_j}R_\mu=\big(q^{-1}a_ja_{2n+4}\big)R_\mu.
\label{eq:TRmu}
\end{gather}
From (\ref{eq:Res}), (\ref{eq:TRmu}) and Corollary \ref{cor:s=n+1Taj},
we obtain the recurrence relation for $I(a_1,a_2,\ldots,a_{2n+3})$ as
\[
T_{a_j}I(a_1,a_2,\ldots,a_{2n+3})=I(a_1,a_2,\ldots,a_{2n+3})
\prod_{\substack{1\le\ell\le 2n+3\\ \ell\ne j}}\frac{1-a_ja_\ell}{1-qa_{2n+4}^{-1}a_\ell^{-1}}.
\]
By repeated use of the above relation, we obtain
\begin{gather*}
I(a_1,a_2,\ldots,a_{2n+3})
 =
\frac{\prod\limits_{k=1}^{2n+3}\big(qa_{2n+4}^{-1}a_k^{-1}\big)_{2N}}{\prod\limits_{1\le i<j\le 2n+3}(a_ia_j)_{2N}}
I\big(q^Na_1,q^Na_2,\ldots,q^Na_{2n+3}\big)
\\
\phantom{I(a_1,a_2,\ldots,a_{2n+3})}{} =
\frac{\prod\limits_{k=1}^{2n+3}\big(qa_{2n+4}^{-1}a_k^{-1}\big)_{\infty}}{\prod\limits_{1\le i<j\le 2n+3}(a_ia_j)_{\infty}}
\lim_{N\to \infty}I\big(q^Na_1,q^Na_2,\ldots,q^Na_{2n+3}\big)
\end{gather*}
and
\begin{gather*}
\lim_{N\to \infty}I\big(q^Na_1,q^Na_2,\ldots,q^Na_{2n+3}\big)
 =
\left(\frac{1}{2\pi\sqrt{-1}}\right)^{\!\!n}\int_{{\mathbb T}^n}
\prod_{i=1}^n
\big(z_i^2\big)_\infty\big(z_i^{-2}\big)_\infty
\\
\qquad {}
\times
\prod_{1\le i<j \le n}(z_iz_j)_\infty\big(z_iz_j^{-1}\big)_\infty\big(z_i^{-1}z_j\big)_\infty\big(z_i^{-1}z_j^{-1}\big)_\infty
\frac{dz_1}{z_1}\wedge\cdots\wedge\frac{dz_n}{z_n}
 =\frac{2^nn!}{(q)_\infty^n}.
\end{gather*}
This completes the proof. \end{proof}

\pdfbookmark[1]{References}{ref}
\LastPageEnding

\end{document}